\documentclass[10pt]{article}
\usepackage{amssymb}
\usepackage{latexsym}
\usepackage{amsmath}

\newtheorem{theorem}{Theorem}[section]
\newtheorem{proposition}[theorem]{Proposition}
\newtheorem{lemma}[theorem]{Lemma}

\newtheorem{corollary}[theorem]{Corollary}

\newtheorem{remark}[theorem]{Remark}

\newcommand{\ENDproof}{\end{proof}$\hspace*{\fill} \square$\\ \smallskip\par\noindent}
\newcommand{\field}{\mathbb K}
\newcommand{\rmv}[1]{}

\newcommand{\dprime}{\prime\prime}
\newcommand{\abs}[1]{\left|#1\right|}

\begin{document}

\title{Fast Jacobian group operations for $C_{3,4}$ curves over a large
  finite field}

\author{Fatima K. Abu Salem\\
 Computer Science Department\\
 American University of Beirut\\
 \texttt{fa21@aub.edu.lb} \\
        and\\
 Kamal Khuri-Makdisi\\
 Mathematics Department and\\
 Center for Advanced Mathematical Sciences\\
 American University of Beirut\\
 \texttt{kmakdisi@aub.edu.lb}
  }
 \date{May 14, 2007}

\maketitle

\begin{abstract}
Let $C$ be an arbitrary smooth algebraic curve of genus $g$ over a large
finite field ${\mathbb K}$.  We revisit fast addition algorithms in the
Jacobian of $C$ due to Khuri-Makdisi (math.NT/0409209, to appear in {\it
Mathematics of Computation}). The algorithms, which reduce to linear
algebra in vector spaces of dimension $O(g)$ once $|{\mathbb K}| \gg g$
and which asymptotically require $O(g^{2.376})$ field operations using
fast linear algebra, are shown to perform efficiently even for certain
low genus curves. Specifically, we provide explicit formulae for
performing the group law on Jacobians of $C_{3,4}$ curves of genus $3$.
We show that, typically, the addition of two distinct elements in the
Jacobian of a $C_{3,4}$ curve requires $117$ multiplications and 2
inversions in ${\mathbb K}$, and an element can be doubled using $129$
multiplications and 2 inversions in ${\mathbb K}$.  This represents an
improvement of approximately 20\% over previous methods.
\end{abstract}

\textbf{Remark (added August 22, 2007):}
A revised version of this article has been published as 
LMS J. Comput. Math. 10 (2007) 307--328
with an appendix of sample \textsc{Magma} code of our algorithms.
The URL for the published version is:

\texttt{http://www.lms.ac.uk/jcm/10/lms2006-049/}

\bigskip

\section{Introduction and background}
\label{introbackground}

This article presents the fastest algorithms to date for arithmetic in the
Jacobians of certain nonhyperelliptic genus 3 curves --- specifically,
$C_{3,4}$ curves over a very large finite field $\field$ that is not of
characteristic $2$ or $3$.  We attain this by adapting the ideas from the
asymptotically fastest algorithms known for general curves of large genus
\cite{KKM,KKM04}.  Those algorithms boil down to linear algebra on
matrices of size $O(g)\times O(g)$, where $g$ is the genus of the curve
(more accurately,
$O\bigl(g(1 + \log g/\log \abs{\field})\bigr) \times O(g)$, but recall that
$\abs{\field}$ is large), and thus have a complexity of $O(g^{2.376})$
using the current record for fast linear algebra.  Our results in this
article illustrate how the asymptotic improvements introduced in
\cite{KKM04}, coupled with further new techniques, actually result in a
significant speedup even for low genus curves that are slightly ``special''
for their genus.  However, fairly special curves, such as hyperelliptic
curves for example, are still probably better implemented using Cantor's
algorithm or the general methods of \cite{Hess}, which have complexity
$O(g^2)$ for curves of bounded gonality, but which have complexity $O(g^4)$
for ``most'' curves of genus $g$.

Previous work on Jacobian group arithmetic for nonhyperelliptic genus $3$
curves includes \cite{BEFG, FOR04}, building on earlier work for curves of
the form $y^3 = x^3 + \alpha x+ \beta$, \cite{BEFG-article, FO-previous}.
The papers \cite{BEFG, FOR04} give slower algorithms for $C_{3,4}$ curves
than ours, under the same hypotheses on $\field$.  This article follows
the lead introduced by \cite{BEFG-article}, and adopted by
\cite{BEFG, FO-previous, FOR04}, in that we present algorithms that are
designed to work only for ``typical'', i.e., sufficiently generic, elements
of the Jacobian of $C$.  Here, non-typical elements belong to a proper
subvariety of the Jacobian, and so occur with frequency
$O(1/\abs{\field})$, which means that they do not arise in practice.  As in
those previous articles, we also measure the complexity of our algorithms
by counting the number of multiplications and inversions that need to be
performed in $\field$.  This is reasonable, because in practical
implementations of finite field arithmetic, addition and subtraction are
much faster than multiplication or inversion, and inversion can take
between $3$ and $10$ times as long as multiplication, as pointed out in
\cite{BEFG}.  Our approach requires $117$ multiplications and $2$
inversions in $\field$ to add a typical pair of distinct elements of the
Jacobian; we abbreviate this complexity as $117M, 2I$.  In contrast, the
complexity of adding a typical pair of distinct elements in \cite{FOR04} is
$145M, 2I$, while the complexity in \cite{BEFG} is $150M, 2I$.  As for
doubling a typical element of the Jacobian, our approach requires $129M,
2I$, as opposed to the doubling algorithm in \cite{FOR04}, which needs
$167M, 2I$, and to that in \cite{BEFG}, which needs $174M, 2I$.  Our
algorithms and those of \cite{FOR04} actually compute first the negative of
a sum of two elements of the Jacobian (respectively $-2$ times an element
during doubling), and then invert the final result.  The final inversion
costs $7M$ in our approach, and $16M$ in \cite{FOR04} (as gathered from an
inspection of their computer code).  This final inversion is not needed if
one wishes to compute a large multiple of an element of the Jacobian by the
usual ``double and add'' method; one can use instead the approach in
\cite{base-2}, which uses the ``addflip'' primitive
$\xi, \xi' \mapsto -(\xi + \xi')$ (where $\xi$ may equal $\xi'$, for
multiplication by $-2$) instead of the usual addition and doubling.  Due to
recent progress in index calculus methods for discrete logarithms (see
\cite{Diem}, \cite{DT06}, and their references), it appears unlikely
that the discrete 
logarithm problem in Jacobians of $C_{3,4}$ curves is worth using as a
cryptographic primitive; the methods of this paper might still be useful
for cover attacks on discrete logarithms of other curves.

For the general problem of computing effectively in Jacobians, our
results in this article confirm the advantages of using the approach of
\cite{KKM, KKM04}.  Even though we write 
down polynomials in this article, our algorithms work mainly via linear
algebra in spaces of sections of line bundles, which we discuss here in the
language of Riemann-Roch spaces $\mathcal{L}(D)$ associated to appropriate
divisors on $C$.  We perform almost no polynomial arithmetic, and instead
use linear algebra on small matrices (essentially, $3 \times 5$ and
$8 \times 10$, both explicitly and implicitly) which are often fairly
structured. For example, our matrix may have two blocks that are almost in
echelon form; hence an intelligent approach to Gaussian elimination
produces efficient algorithms.  We also optimise by hand any parts of the
calculations that yield easily to an ad hoc trick or to more
systematic approaches. We hope that some of these methods can be useful
elsewhere.

The second named author (KKM) would like to thank K. and C. Adal for
providing computer access and an agreeable work environment during the
summer of 2006, when the author's usual office was inaccessible.

\section{Overview of our algorithms}
\label{overview}

Consider a $C_{3,4}$ curve $C$ of genus $3$ over a large finite field
$\field$ with $q = p^n$ elements. We assume that $p$, the characteristic of
$\field$, is neither $2$ nor $3$ (similarly to \cite{BEFG, FOR04}; those
articles also exclude characteristic $5$). Let $P_{\infty} \in C$ denote
the distinguished point at infinity and $D$ a $\field$-rational divisor on
$C$. Write $\mathcal{L}(D)$ for the Riemann-Roch space of rational
functions on $C$ with prescribed zeros and poles at $D$:
\[
\mathcal{L}(D) = \{ F \in \field(C) \left \vert \right. (F) \geq -D\}.
\]
Write $\mathcal{R}$ for the affine coordinate ring of $C-\{P_{\infty}\}$;
hence $\mathcal{R} = \cup_{N \geq 0} \mathcal{L}(NP_{\infty})$.
By definition of a $C_{3,4}$ curve, $\mathcal{R}$ is generated as a
$\field$-algebra by two elements $x$, $y$ whose valuations $v_{P_{\infty}}$
are given by 
\[
\begin{array}{rcl}
v_{P_{\infty}}(x) & = & -3\\
v_{P_{\infty}}(y) & = & -4.\\
\end{array}
\]
The only relation between $x$ and $y$ is a $\field$-linear dependence
$f(x,y) = 0$ between $1$, $x$, $y$, $x^2$, $xy$, $y^2$, $x^3$, $x^2y$,
$xy^2$, $y^3$, $x^4 \in \mathcal{L}(12P_{\infty})$. Thus, the affine
coordinate ring of $C-\{P_{\infty}\}$ is $\mathcal{R} =
\field[x,y]/(f(x,y))$. After a change of variables of the form
\[
\left \{
\begin{array}{rcl}
x & \mapsto & u_1x + u_2\\
y & \mapsto & u_3y + u_4x + u_5,
        \quad u_1,\ldots,u_5 \in \field, \quad u_1,u_3 \neq 0,
\end{array}
\right.
\]
we can assume that the equation of the curve is
\begin{equation}
f(x,y) = y^3 - x^4 + p_2x^2y + p_1xy + p_0 y + q_2 x^2 + q_1x + q_0 = 0.
\label{eqn}
\end{equation}
We further write $W^{N} = \mathcal{L}(NP_{\infty})$; it is the subspace of
$\mathcal{R}$ spanned by the monomials
\[
\{ x^iy^j \left \vert \right. 3i + 4j \leq N\},
\]
subject to the relation \eqref{eqn}.  To obtain a basis of $W^N$, we
restrict ourselves to monomials with exponent pairs $(i,j)$ with 
$j \leq 2$, or alternatively to pairs $(i,j)$ with $i \leq 3$; this takes
equation (\ref{eqn}) into account. Note that
\[
\begin{array}{rcl}
W^0 & = & W^1\\
&  = & W^2 =  \field \cdot 1 \mbox{ is }1\mbox{-dimensional}\\
W^3 & = & \field \cdot 1 + \field \cdot x \mbox{ is }2\mbox{-dimensional}\\
W^4 & = & W^5 \\
& = & \field \cdot 1 + \field \cdot x + \field \cdot y
                                           \mbox{ is }3\mbox{-dimensional}\\
\end{array}
\]
and for $N \geq 6$, $W^N$ is $(N-2)$-dimensional.

Let $D$ be an effective $\field$-rational divisor.   Following the approach
of \cite{KKM, KKM04}, we represent $D$ by the space $W^{N}_{D}$ defined by
\[
W^{N}_{D} = \mathcal{L}(NP_{\infty} - D) \subset W^N
\]
for some suitable positive integer $N$. If $D$ is arbitrary of degree $d$,
then we need to consider $N \geq d + 6$ (here, $6 = 2g$ for $g = 3$, the
genus of the curve, to ensure that $W^{N}_{D}$ is base-point
free). However, for a typical divisor $D$, we can take $N = d+4$ (here,
$4 = g+1$). This is a consequence of the following standard result from the
theory of linear series on curves: 
\begin{proposition}
Let $D$ be a typical effective $\field$-rational divisor of degree
$d \geq 3$ on $C$. In particular, $P_{\infty}$ does not belong to the 
support of $D$. Then 
\[
\dim W^{N}_{D} = \left \{
\begin{array}{rcl}
    0  & \mbox{{\rm if }} & N \leq d+2\\
 N-d-2 & \mbox{{\rm if }} & N \geq d+2.
\end{array} 
\right.
\]
Furthermore, if $N \geq d+4$, then $W^{N}_{D}$ is base-point free, and
there exist two elements $F \in W^{d+3}_{D}$ and
$G \in W^{d+4}_{D} - W^{d+3}_{D}$ that form a basis for the $2$-dimensional
subspace $W^{d+4}_{D} \subset W^{N}_{D}$, with the property that the only
common vanishing of $F$ and $G$ occurs at $D$.  In other words,
\[
\begin{array}{rcl}
(F) & = & -(d+3)P_{\infty} + D + E\\
(G) & = & -(d+4)P_{\infty} + D + E^{\prime},
\end{array}
\]
where $E$ and $E^{\prime}$ are disjoint effective divisors.
\end{proposition}
\begin{remark}
Since $F$ and $G$ above vanish simultaneously only at $D$, we see that our
basis $\{F,G\}$ for $W^{d+4}_{D}$ is in fact an ideal generating set (an
\emph{IGS}) for $D$ in the terminology of \cite{KKM04}. Thus, the ideal
$\langle F,G \rangle = \mathcal{R}F + \mathcal{R}G$ of the affine
coordinate ring $\mathcal{R}$ is the ideal of elements of $\mathcal{R}$
(i.e., of regular functions on $C - \{P_\infty\}$) vanishing on $D$. The
quotient $\mathcal{A} = \mathcal{R}/\langle F,G \rangle$ is a
$d$-dimensional $\field$-algebra describing the ``values'' that a
polynomial can take at the points of $D$. This makes sense even if the
points of $D$ are not all defined over $\field$, so long as the divisor $D$
itself is $\field$-rational. Moreover, there is a $\field$-linear map 
\[
W^{N}/W^{N}_{D} \to \mathcal{A}
\]
that is a bijection for $N \geq d+2$, for typical $D$ with $d \geq 3$.
\label{rem1}
\end{remark}
\begin{remark}
As mentioned above, a ``typical'' divisor $D$ is one that does not belong
to a specific proper (that is, at most $(d-1)$-dimensional) subvariety of
the $d$-dimensional symmetric power ${\rm Sym}^{d}C$ parametrising the
degree $d$ effective divisors on $C$. For very large $q = \abs{\field}$,
the probability for a divisor $D$ to be non-typical is $O(1/q)$.  For
enormous $q$, we do not expect to ever chance upon a non-typical divisor in
our calculations. In case we do, it was already remarked in
\cite{BEFG, BEFG-article} that we can then use a slower algorithm that
works for all divisors. For example, we can use the larger space
$W^{d+6}_{D}$ instead of $W^{d+4}_{D}$, and adapt the algorithms
accordingly.
\end{remark}

We now discuss how we compute with typical elements of the Jacobian $J$ of
$C$. An element $\xi \in J(\field)$ can be represented as the divisor class
$[D-3P_{\infty}]$ for some effective $\field$-rational divisor $D$ with
$\deg D = 3$. A typical class corresponds to a typical divisor $D$ in a
unique way. In turn, we represent $D$ by a basis $\{F,G\}$ for the
$2$-dimensional space $W^{7}_{D}$; i.e., by elements $F,G \in W^{7}_{D}
\subset \mathcal{R} = \field[x,y]/(f(x,y))$. We can choose the basis
$\{F,G\}$ to have the form 
\begin{equation}
\left \{
\begin{array}{rclcl}
F & = & x^2 + ay + bx + c & \in & W^{6}_{D} \subset W^{7}_{D}\\
G & = & xy  + dy + ex + f & \in & W^{7}_{D} - W^{6}_{D}.
\end{array}
\right.
\label{star0}
\end{equation}
Here $a \neq 0$ for typical divisors, and, for technical reasons, we also
store the inverse $a^{-1}$ along with the coefficients $a,b,\ldots,f \in
\field$ in order to represent $\xi = [D-3P_\infty]$.

Our addition algorithm begins with a typical pair $\xi,\xi^{\prime} \in
J(\field)$ and computes their sum $\xi+\xi'$.  Our doubling algorithm
corresponds to the special case $\xi = \xi'$, in which case we compute
$2\xi = \xi + \xi'$.  In both cases, we first compute $\xi^{\dprime} =
-(\xi+\xi^{\prime})$, the ``addflip'' of the two elements in the
terminology of \cite{KKM, KKM04}.  We then compute $\xi''' = -\xi''$.  In
practice, most of the use of Jacobian arithmetic will be to find a multiple
$m \cdot \xi$ with $m \in \mathbb{Z}$.  In that case, we can use the ``base
$-2$ expansion'' of \cite{base-2} and only find the addflips $\xi''$ in the
intermediate steps without any need for further negations. 

We thus start with $\xi = [D-3P_{\infty}]$ and $\xi^{\prime} =
[D^{\prime}-3P_{\infty}]$, with bases $\{F,G\}$ for $W^{7}_{D}$ and
$\{F^{\prime},G^{\prime}\}$ for $W^{7}_{D^{\prime}}$. In our first phase
(Steps $1$ and $2$ below) we produce a basis $\{F^{\dprime},G^{\dprime}\}$
for $W^{7}_{D^{\dprime}}$ where $[D+D^{\prime}+D^{\dprime}-9P_{\infty}] =
0$ in $J(\field)$. Thus $F'', G''$ represent $\xi^{\dprime} =
[D^{\dprime}-3P_{\infty}] = -(\xi+\xi^{\prime})$.  In our second phase
(Step $3$ below), we find a basis $\{F''', G'''\}$ for $W^7_{D'''}$, where
$[D''+D'''-6P_\infty]=0$ in $J(\field)$. At this point, $F''',G'''$
represent $\xi''' = [D'''-3P_{\infty}] = -\xi''$.  Along the way, we also
obtain the inverses $(a'')^{-1}$ and $(a''')^{-1}$ of the analogous
coefficients in $F''$ and $F'''$. Here is a more detailed overview:

\subsection{Step 1}
\label{step1}
This step comprises Sections \ref{prelim}--\ref{findST} of this article. We
first determine the space $W^{10}_{D+D^{\prime}}$ along with its subspace
$W^{9}_{D+D^{\prime}}$. Since $D+D^{\prime}$ is typical, we have that
$\dim W^{9}_{D+D^{\prime}} = 1$ and $\dim W^{10}_{D+D^{\prime}} = 2$. Thus,
there exists a basis $\{s,t\}$ for $W^{10}_{D+D^{\prime}}$ of the form
\begin{equation}
\begin{array}{rcl}
s & = & x^3 + s_1y^2 + s_2xy + s_3x^2 + s_4y + s_5x + s_6\\
& = & 0x^2y + 1x^3 + \ldots \quad \in W^{9}_{D+D^{\prime}}
                                      \subset  W^{10}_{D+D^{\prime}}\\
t & = & x^2y + t_1y^2 + t_2xy + t_3x^2 + t_4y + t_5x + t_6\\
& = & 1x^2y + 0x^3 + \ldots  \quad
                     \in W^{10}_{D+D^{\prime}} -  W^{9}_{D+D^{\prime}},
\end{array}
\label{dstar0}
\end{equation}
with $s_1,\ldots,s_6, t_1,\ldots,t_6 \in \field$. Our aim is thus to find
$s$ and $t$. Note that the principal divisor $(s)$ has the form $(s) =
D + D^{\prime} + D^{\prime\prime} - 9P_{\infty}$ for some effective
$\field$-rational divisor $D''$ of degree $3$. Hence, 
$[D + D' + D'' - 9P_{\infty}] = 0$, and $\xi'' = -(\xi+\xi')$ as desired.

Carrying out Step $1$ depends on whether $D \neq D^{\prime}$ (corresponding
to addition) or $D = D^{\prime}$ (corresponding to doubling).

\subsubsection{Point addition}
If $D \neq D^{\prime}$, then $D$ and $D^{\prime}$ typically have no point
in common, in which case
\[
W^{10}_{D+D^{\prime}} = W^{10}_{D} \cap W^{10}_{D^{\prime}}.
\]
We find this intersection by looking for those elements of
$W^{10}_{D^{\prime}}$ that map to zero in the quotient ring
$\mathcal{A} = \mathcal{R}/\langle F,G \rangle$ (hence such elements also
vanish at $D$). We set up $\mathcal{A}$ in Section \ref{prelim}, compute
how a basis for $W^{10}_{D^{\prime}}$ maps to $\mathcal{A}$ in 
Section \ref{pre-add}, and find the kernel of the map 
($W^{10}_{D^{\prime}} \to \mathcal{A}$) in Sections \ref{findkernel} and
\ref{findST}.  

\subsubsection{Point doubling}
If $D = D^{\prime}$, then we compute $W^{10}_{2D}$ as the subspace of
elements $L \in W^{10}_{D}$ whose differential $dL$ also vanishes at
$D$. This differs from the case of addition above only in computing a map
$(W^{10}_{D} \to \mathcal{A}^{\prime}): L \mapsto dL$ ``mod''
$\langle F,G \rangle$, where $\mathcal{A}^{\prime}$ is a $3$-dimensional
$\field$-vector space describing the ``values'' that $dL$ can take at the
points of $D$.  We describe this in Section \ref{pre-double}, the analogue
of Section \ref{pre-add} with respect to point addition.  Thereafter, the
remaining calculations in Sections \ref{findkernel} and \ref{findST} proceed
similarly to the case of point addition.

\subsection{Step 2}
\label{step2}

This step comprises Sections \ref{xyst} and \ref{last} below. At this
stage, we have a basis $\{s,t\}$ for $W^{10}_{D+D^{\prime}}$ as in
(\ref{dstar0}), which is typically an IGS for $D+D^{\prime}$ as in Remark
\ref{rem1}. Thus,
\[
\begin{array}{rcl}
(s) & = & D + D^{\prime} + D^{\dprime} - 9P_{\infty},\\
(t) & = & D + D^{\prime} + E^{\dprime} - 10P_{\infty},
\end{array}
\]
with $D^{\dprime}$ and $E^{\dprime}$ disjoint. We note that $sW^8 =
W^{17}_{D+D^{\prime}+D^{\dprime}}$ as in \cite{KKM}.  Taking a basis of
monomials for $W^{8}$, we see that the following is a basis for $sW^{8}$: 
\[
\{ s, xs, ys, x^2s, xys, y^2s\}.
\]
We next compute $W^{7}_{D^{\dprime}}$. It is the ``quotient'', as in
\cite{KKM04}, of $sW^8 = W^{17}_{D+D^{\prime}+D^{\dprime}}$ by the IGS
$\{s,t\}$ for $D+D^{\prime}$: 
\begin{equation}
\begin{array}{rcl}
W^{7}_{D^{\dprime}} & = & sW^8 \div \{s,t\}\\
& = & \{\ell \in W^7 \left \vert \right. s\ell, t\ell \in sW^8 \}\\
& = & \{ \ell \in W^7 \left \vert \right. t\ell \in sW^8\}.
\end{array}
\label{tstar0}
\end{equation}
Since $W^7$ has basis $\{1,x,y,x^2,xy\}$ and we have a basis for $sW^8$,
the condition $t\ell \in sW^8$ amounts to finding a linear combination of
$t$, $xt$, $yt$, $x^2t$, and $xyt$ that is also a linear combination of
$s$, $xs$, $ys$, $x^2s$, $xys$, and $y^2s$. Equivalently, we must determine
the intersection of the $5$- and $6$-dimensional subspaces $tW^7$ and
$sW^8$ inside $W^{17}$. This intersection will have a basis of the form
$\{tF^{\dprime},tG^{\dprime}\}$, where $\{F^{\dprime},G^{\dprime}\}$ are a
basis for the space $W^{7}_{D^{\dprime}}$ of solutions for $\ell$ in
(\ref{tstar0}) above. Note that the intersection appears to take place in
the $15$-dimensional space $W^{17}$ (where typical $5$ and $6$-dimensional
spaces do not intersect), but actually occurs inside the $9$-dimensional
space $W^{17}_{D+D^{\prime}}$, which contains (in fact, is generated by)
the two subspaces $tW^7$ and $sW^8$. This reduces the amount of linear
algebra that we need to perform. We formalise this in the following lemma: 
\begin{lemma}
\label{lemma2.1}
Let $\ell \in W^7$. Then $t\ell \in sW^8$ if and only if $t\ell \in
sW^8 + W^9$. (This is equivalent to saying that $t\ell$ is congruent to an
element of $sW^8$ in the quotient space $W^{17}/W^{9}$.) 
\end{lemma}

\begin{proof}
Trivially, $t\ell \in sW^8$ implies that $t\ell \in sW^8 + W^9$. To prove
the converse, suppose that $t\ell = s\ell^{\prime} + \ell^{\dprime}$,
with $\ell^{\prime} \in W^8$ and $\ell^{\dprime} \in W^9$. Note that
$t\ell, s\ell^{\prime} \in W^{17}_{D+D^{\prime}}$. Then, since
$\ell^{\dprime} \in W^9$, we obtain
\[
\ell^{\dprime} = t\ell - s\ell^{\prime} \in W^{9}_{D+D^{\prime}}
               = \field \cdot s,
\]
and so we can write
\[
t\ell - s\ell^{\prime} = \alpha s, \quad \alpha \in \field,
\]
from which we have
\[
t\ell = (\ell^{\prime}+\alpha)s \in sW^8,
\]
as required. Note incidentally that $sW^8 \cap W^9 = \field s$, so
$\dim (sW^8 + W^9) = 6 + 7 - 1 = 12$. 
\ENDproof 
We conclude from the above discussion that we can obtain $F^{\dprime},
G^{\dprime} \in W^{7}_{D^{\dprime}}$ as follows:
\begin{enumerate}
\item
Denote $F^{\dprime}$ or $G^{\dprime}$ by
$\ell = d_1+d_2x + d_3y + d_4x^2 + d_5xy$.  Here $\{d_4, d_5\} = \{0,1\}$
in some order, and we must solve for $d_1, d_2, d_3$ such that 
$t\ell \in sW^8 + W^9$. 
\item
Find $\overline{C_1},\ldots,\overline{C_5}$, the images of
$t,xt,yt,x^2t,xyt$ in the $3$-dimensional quotient space
$W^{17}/(sW^8+W^9)$. (One can moreover see from Section \ref{last} that a
basis for this quotient space is given by the images of $x^2y,xy^2$, and
$x^2y^2$.)
\item
The three resulting equations
$d_1\overline{C_1} + \ldots + d_5\overline{C_5} = 0$ allow us (in the
typical case) to express $d_1,d_2,d_3$ in terms of $d_4,d_5$. We thus get a
basis $\{(c^{\dprime},b^{\dprime},a^{\dprime},1,0),
         (f^{\dprime},e^{\dprime},d^{\dprime},0,1)\}$
for the space $\{(d_1,\ldots,d_5) \mid
d_1\overline{C_1} + \ldots + d_5\overline{C_5} = 0\}$. This corresponds
to elements $F^{\dprime} = c^{\dprime}+b^{\dprime}x + a^{\dprime}y+x^2$ and
$G^{\dprime} = f^{\dprime} + e^{\dprime}x + d^{\dprime}y + xy$ that form a
basis for $W^{7}_{D^{\dprime}}$.  The structure of the system of linear
equations allows us to find $(a'')^{-1}$ along the way at minimal extra
cost.
\end{enumerate}

\subsection{Step 3}
\label{step3}
This step comprises Section \ref{negation}.  At this point we have obtained
our IGS $\{F'', G''\}$ for the divisor $D''$, where
$\xi'' = [D''-3P_{\infty}] = - (\xi + \xi')$.  We also know $(a'')^{-1}$.
We now discuss how to negate this to obtain $\xi''' = -\xi'' = \xi + \xi'$.
The divisor of $F''$ has the form $(F'') = D'' + D''' - 6P_{\infty}$ for
some effective $\field$-rational divisor $D'''$, and it follows that
$\xi''' = [D'''-3P_{\infty}]$.  We thus seek the polynomials 
 \[
\begin{array}{rcl}
F''' & = & x^2 + a'''y + b''' x + c''' \in W^6_{D'''}\\
G''' & = & xy + d'''y + e''' x + f''' \in W^7_{D'''}
\end{array}
\]
that represent $D'''$ and hence $\xi'''$.  We easily observe that
$F'' = F'''$, since $W^6_{D'''} = W^6_{D''} = W^6_{D''+D'''} = 
\field \cdot F''$.  Hence $a''' = a''$, so we trivially know the inverse
$(a''')^{-1}$.

It remains to find $G'''$.  Analogously to \eqref{tstar0} and to
Lemma \ref{lemma2.1}, we have $F''W^8 = W^{14}_{D'' + D'''}$, and so  
\begin{equation}
\begin{array}{rcl}
W^{7}_{D'''} & = & F''W^8 \div \{F'',G''\}\\
& = & \{\ell \in W^7 \left \vert \right. G''\ell \in F''W^8 \}\\
& = & \{ \ell \in W^7 \left \vert \right. G''\ell \in F''W^8 + W^6\}.
\end{array}
\label{tstar1}
\end{equation}
We thus have $G''G'''+F''H=0$ for some $H \in W^8$. We can in principle
carry out an analogous computation to Step 2, but this case is small enough
that it is worth our while to carry out the calculation directly and to
hand-optimise it to find $G'''$.  We also find an explicit expression for
$H$, which is useful in a different context that we encounter in Section
\ref{pre-double}.

\section{Preliminary to both point addition and doubling}
\label{prelim}

Consider the input $F = x^2 + ay + bx + c$ and $G = xy + dy + ex + f \in
W^{7}_{D}$ representing a divisor $D$ of degree $3$.  We know that 
$\langle F,G \rangle = \mathcal{R}F + \mathcal{R}G$ is the ideal of regular
functions on $C - \{P_\infty\}$ vanishing at $D$. Our goal is to be able to
compute in the algebra of ``values'' of polynomials at $D$, given by 
\[
\mathcal{A} = \mathcal{R} / \langle  F,G \rangle.
\]
Since $\deg D = 3$, we have $\dim_{\field} \mathcal{A} = 3$. Given
$u \in \mathcal{R}$, the element $\overline{u} \in \mathcal{A}$ denotes the
reduction of $u$ modulo $\langle F,G \rangle$. 

\begin{lemma}
A $\field$-basis for $\mathcal{A}$ is
$\{\overline{1},\overline{x},\overline{y}\}$. Furthermore,
\begin{equation}
\overline{x}^2 = -a \overline{y} - b \overline{x} - c \overline{1}
\label{s1}
\end{equation}
\begin{equation}
\overline{x}\overline{y}
              = - d \overline{y} - e \overline{x} - f \overline{1}
\label{s2}
\end{equation}
\begin{equation}
\overline{y}^2 = -g \overline{y} - h \overline{x} - i \overline{1}
 \label{s3}
\end{equation}
where $a,b,c,d,e$ and $f$ are the coefficients of $F$ and $G$ 
and 
\begin{equation}
\label{s4}
\begin{array}{rcl}
g & = &a^{-1}\left (c + d(d-b)\right) + e\\
h & = &a^{-1}(ed -f)\\
i & = & a^{-1}\left (ec +f(d-b)\right).
\end{array}
\end{equation}
\end{lemma}
\begin{proof}
Equations (\ref{s1}) and (\ref{s2}) reflect the fact that $F,G \in
\langle F,G \rangle$. Equations (\ref{s3}, \ref{s4}) come from expanding 
$(y+e)F - (x+b-d)G \in \langle F,G \rangle$. Equations
(\ref{s1}, \ref{s2}, \ref{s3}) show that every element 
$\overline{u} \in \mathcal{A}$ can be written as a $\field$-linear
combination of $\overline{1},\overline{x}$ and $\overline{y}$. Since
$\mathcal{A}$ is three dimensional, we obtain that $\overline{1}$,
$\overline{x}$ and $\overline{y}$ are linearly independent.
\ENDproof

Given $u \in \mathcal{R}$, we represent its reduction $\overline{u} =
\alpha \overline{1} + \beta \overline{x}+ \gamma \overline{y} \in
\mathcal{A}$ by the column vector 
\[
B_{u} = \left (\begin{array}{l} \alpha \\ \beta \\ \gamma \end{array} 
        \right)
   \in \field^{3}.
\]
We then have
\begin{proposition}
Assume given $F$ and $G$, as well as the inverse $a^{-1}$.
\begin{enumerate}
\item 
For $B_{u}$ defined as above, we have
\[
B_{xu} = T_xB_u,
\qquad\qquad
B_{yu} = T_yB_u,
\]
where $T_x$ and $T_y$ are the matrices of multiplication by $x$ and $y$
on $\mathcal{A}$, with respect to the ordered basis
$\{\overline{1}, \overline{x}, \overline{y}\}$: 
\[
T_x = \left (
\begin{array}{lll}
0 & -c & -f\\
1 & -b & -e\\
0 & -a & -d
\end{array}
\right),
\qquad\qquad
T_y = \left (
\begin{array}{lll}
0 & -f & -i\\
0 & -e & -h\\
1 & -d & -g
\end{array}
\right).
\]
\item
We have the entries of $T_{x}$ for free (i.e., at a cost of $0M$);
multiplying $T_x$ by a vector $B_u$ costs $6M$. 
\item
We can compute the entries of $T_{y}$ using $7M$.  Once we know $T_y$,
multiplying $T_y \cdot B_u$ to get $B_{yu}$ also costs $6M$. 
\item
If we do not already know $T_y$, we can obtain $B_{yu}$ directly at a cost
of $11M$. 
\end{enumerate}
\label{cost1}
\end{proposition}
\begin{proof}
The proof of parts 1--3 is immediate by inspecting (\ref{s1}--\ref{s4})
above.  As for part 4, we need to compute the reduction modulo $\langle
F,G\rangle$  of $v = \alpha y + \beta xy + \gamma y^2$ in order to obtain
$B_{yu}$.  Now $v$ is congruent to $w = v - \gamma a^{-1} (yF - xG)$, so we
have 
\[
w = \gamma a^{-1} f x + (\alpha - \gamma a^{-1} c) y
  + \gamma a^{-1} e x^2 + [\beta - \gamma a^{-1}(b-d)] xy
  = \delta x + \varepsilon y + \zeta x^2 + \eta xy
\]
where $\delta, \varepsilon, \zeta, \eta$ can be calculated using $5M$
(first find $\gamma a^{-1}$).  Then the reduction modulo $\langle
F,G\rangle$ of $v$ is $w - \zeta F - \eta G$, whence 
\[
B_{yu} =
\begin{pmatrix}
0 \\ \delta \\ \varepsilon
\end{pmatrix}
-
\zeta
\begin{pmatrix}
c \\ b \\ a
\end{pmatrix}
-
\eta
\begin{pmatrix}
f \\ e \\ d
\end{pmatrix},
\]
costing an additional $6M$.
\ENDproof

\section{First stage of addition of two distinct divisor classes: setting
up a system of equations whose solution will determine
$W^{10}_{D+D^{\prime}}$}
\label{pre-add}

Our input is now the descriptions of two typical degree $3$ divisors $D,
D'$, given by $F = x^2 + ay + bx + c, G = xy + dy + ex + f \in W^{7}_{D}$
and $F^{\prime} = x^2 + a^{\prime}y + b^{\prime}x + c^{\prime},
G^{\prime} = xy + d^{\prime}y + e^{\prime}x + f^{\prime} \in
W^{7}_{D^{\prime}}$, along with the inverses $a^{-1}$ and $(a')^{-1}$.  Our
goal in this section is to determine a $3 \times 5$ matrix $M$ whose five
columns are respectively $B_{F^{\prime}}, B_{xF^{\prime}}, B_{yF^{\prime}},
B_{G^{\prime}}, B_{xG^{\prime}}$, in the notation of Section \ref{prelim}.
The kernel of $M$ will then correspond to $W^{10}_{D+D^{\prime}}$ as
follows: if $v = (c_1, c_2, c_3, c_4, c_5)^\mathbf{T}$ is a (column) vector
in $\field^5$, then we identify it with the linear combination
\[
L = (c_1+c_2x+c_3y)F^{\prime} + (c_4+c_5x)G^{\prime}
   \in \langle F',G' \rangle  \cap W^{10} = W^{10}_{D'}.
\] 
Then $M v = 0$ if and only if $\overline{L} = \overline{0}$ in
$\mathcal{A}$, which is equivalent to  
\[
 L \in \langle F,G \rangle \cap W^{10}_{D^{\prime}} = W^{10}_{D+D^{\prime}},
\]
where the last equality follows from the fact that $D$ and $D^{\prime}$ are
disjoint.

\begin{proposition}
Given $F,G,F',G', a^{-1}$ as above, we can compute the matrix $M$ at a cost
of $22M$. 
\label{add-cost2}
\end{proposition}
\begin{proof}
The first column  $B_{F^{\prime}}$ of $M$ comes from
\[
\begin{array}{rcl}
\overline{F^{\prime}} & \equiv & \overline{F^{\prime}-F}
                                     \bmod \langle F,G \rangle \\
& = & (a^{\prime}-a)y + (b^{\prime}-b) x + (c^{\prime}-c).
\end{array}
\]
Hence we get the following result for free (i.e., $0M$):
\[
B_{F^{\prime}} = \left ( 
\begin{array}{l}
c^{\prime}-c\\
b^{\prime}-b\\
a^{\prime}-a
\end{array}
\right).
\]
We similarly obtain the fourth column $B_{G^{\prime}}$ of $M$ for free:
\[
B_{G^{\prime}} = \left ( 
\begin{array}{l}
f^{\prime}-f\\
e^{\prime}-e\\
d^{\prime}-d
\end{array}
\right).
\]
We now compute the second and fifth columns $B_{xF^{\prime}}$ and
$B_{xG^{\prime}}$ by noting the block matrix equation involving the matrix
$T_x$ of Proposition \ref{cost1}: 
\[
\left ( B_{xF^{\prime}} \mid B_{xG^{\prime}} \right)
  = T_{x} \left ( B_{F^{\prime}}\mid B_{G^{\prime}} \right). 
\]
Since the first column of $T_{x}$ is $\left (0,1,0\right)^{\bf{T}}$, its
interaction with the first row of 
$\left ( B_{F^{\prime}} \mid B_{G^{\prime}} \right)$ can be computed
without any multiplication in $\field$. We must then multiply the $3 \times
2$ submatrix consisting of the second and third columns of $T_{x}$ with the
$2 \times 2$ submatrix consisting of the second and third rows of
$\left (B_{F^{\prime}} \mid B_{G^{\prime}} \right)$. This can be
done using $11M$ using a Strassen's type multiplication on a $2 \times 2$
sub-block, which saves one multiplication over the ``naive'' method.
Finally, we use part 4 of Proposition \ref{cost1} to compute the third
column $B_{yF^{\prime}}$ from $B_{F^{\prime}}$ at a further cost of $11M$.
This concludes the proof.
\ENDproof

\section{First stage of doubling a divisor class: setting up a system of
equations whose solution will determine $W^{10}_{2D}$} 
\label{pre-double}

In this section, we take $D' = D$, so our input consists of the two
polynomials $F = x^2 + ay + bx + c, G = xy + dy + ex + f \in W^{7}_{D}$,
where $D$ is a degree $3$ divisor.  Analogously to Section \ref{pre-add},
we will construct a $3 \times 5$ matrix, which we also label as $M$, whose
columns represent the ``reductions modulo $\langle F,G \rangle$'' of the
differential forms $dF, d(xF), d(yF), dG, d(xG)$.  These differential forms
are regular on $C - \{P_\infty\}$, so we really want the columns of $M$ to
represent the ``values'' of $dF, \dots, d(xG)$ at the points of $D$, much
in the same way that elements of the algebra $\mathcal{A}$ describe values
at $D$.

As in Section \ref{pre-add}, a column vector
$v = (c_1, c_2, c_3, c_4, c_5)^\mathbf{T} \in \field^5$ represents 
\[
L = (c_1 + c_2 x + c_3 y)F + (c_4 + c_5 x)G
         \in \langle F,G \rangle \cap W^{10} = W^{10}_{D}.
\]
This time, $Mv = 0$ if and only if the differential form
$dL = c_1 \,dF + c_2 \,d(xF) + c_3 \,d(yF) + c_4 \,dG + c_5 \,d(xG)$
vanishes at $D$.  Since generically the points of $D$ are distinct, this
means that such an $L$ vanishes to second order at the points of $D$, so we
obtain that $Mv = 0$ if and only if $L \in W^{10}_{2D}$.  Since, e.g.,
$d(xF) = x\,dF + F\,dx$, and $F$ vanishes at $D$, we see that the value of
$d(xF)$ at $D$ is the same as that of $x\, dF$, and so forth.  Thus the
columns of our matrix $M$ can be taken to represent suitable ``reductions
modulo $\langle F,G \rangle$''
\[
\overline{dF}, \overline{x\,dF}, \overline{y\,dF},
     \overline{dG}, \overline{x\,dG}
\]
which we need to explain.  We write $d\mathcal{R}$ for the
$\mathcal{R}$-module of differential forms on $C-\{P_{\infty}\}$; then
$d\mathcal{R}$ is generated by $dx$ and $dy$, with the sole relation
$df = 0$ for $f(x,y)$ the equation of the curve in (\ref{eqn}). 
\begin{lemma}
The $\mathcal{R}$-module $d\mathcal{R}$ is free of rank $1$, and is
generated by a differential form $\omega_0$ such that  
\begin{equation}
dx = f_y\omega_0, \qquad dy = -f_x \omega_0, 
\label{tstar}
\end{equation}
where $f_y  =  \partial f/\partial y$ and $f_x =  \partial f/\partial x$.
\end{lemma}
\begin{proof}
The relation $df = 0$ means that 
\begin{equation}
\label{ostar}
f_x \,dx + f_y \,dy = 0.
\end{equation}
Since $C$ is nonsingular, $f,f_x$, and $f_y$ have no common zeros in the
algebraic closure $\overline{\field}$.  We can therefore write 
\begin{equation}
1 = r_1 f_x + r_2 f_y \quad \mbox{for some } r_1,r_2 \in \mathcal{R},
\label{dstar}
\end{equation}
and we define 
\[
\omega_0 = r_2 \, dx - r_1 \, dy \in d\mathcal{R}.
\]
Some algebra with (\ref{ostar}, \ref{dstar}) then implies equation
(\ref{tstar}). In particular, $dx,dy \in \mathcal{R}\omega_0$ so that
$\omega_0$ generates $d\mathcal{R}$ as an $\mathcal{R}$-module. To see that
the annihilator of $\omega_0$ is $0$, one can argue directly from
(\ref{ostar}), (\ref{dstar}) and the definition of $\omega_0$, or one can
use the fact that $d\mathcal{R}$ is a rank one projective module over the
Dedekind domain $\mathcal{R}$, and hence free, as it has a global generator
$\omega_0$.
\ENDproof

At this stage, we can state precisely what we mean by the reduction
modulo $\langle F, G \rangle$ of the differential forms $dF, \dots, x \,dG$.

\begin{corollary}
Define the reduction of an element of $d\mathcal{R}$ to be its image in
$\mathcal{A}' = d\mathcal{R}/\langle F, G \rangle d\mathcal{R}$.  Then
$\mathcal{A}'$ is a free $\mathcal{A}$-module of rank $1$, generated by the
reduction $\overline{\omega_0}$.
\end{corollary}
We can in fact choose any generator $\overline{\omega}$ of $\mathcal{A}'$,
not just $\overline{\omega_0}$.  Then an element of $\mathcal{A}'$ has the
form $\overline{r\omega}$ for some $r \in \mathcal{R}$, where the reduction
$\overline{r} \in \mathcal{A}$ is well-defined.  We then represent a
reduction $\overline{r\omega}$ by the vector $B_r \in \field^3$.  Our
choice of $\overline{\omega}$ below was inspired by a careful reading of
the formulae for doubling in \cite{FOR04}.  This saves us several
multiplications over using the generator $\overline{\omega_0}$. 

\begin{lemma}
For a typical divisor $D$:
\begin{enumerate}
\item
The reduction $\overline{dF}$ generates the $\mathcal{A}$-module
$\mathcal{A}'$. 
 \item
There exist $G_1 \in W^{7}, H_1 \in W^{8}$ such that $FH_1 + GG_1 = 0$, and
$\overline{G_1}$ is a unit in the ring $\mathcal{A}$. 
\item
There exists a generator $\overline{\omega} \in \mathcal{A}'$ such that 
\begin{equation}
\label{dFdG}
\overline{dF} = \overline{G_1\omega}, \qquad
\overline{dG} = -\overline{H_1\omega}. 
\end{equation}
\end{enumerate}
\end{lemma}
\begin{proof}
The first assertion holds because $F$ typically vanishes to order exactly
one at each point of $D$, so $dF$ is nonzero at the points of $D$.  The
second assertion comes from our results in Subsection \ref{step3} and
Section \ref{negation} (replace $\{F'', G'',G''',H\}$ there by 
$\{F, G, G_1,H_1\}$; no circular reasoning is involved).  The divisor of
$F$ is $(F) = D + D_1 - 6P_\infty$ for a ``complementary'' divisor $D_1$ of
$D$, which is typically disjoint from $D$. (In the original setting of
Section \ref{negation}, $D'''$ was the complementary divisor of $D''$).
Moreover, the only points where $F$ and $G_1$ simultaneously vanish are
typically those of $D_1$, since $\{F,G_1\}$ are an IGS for $D_1$ (indeed,
they are a basis for $W^7_{D_1}$).  Thus $G_1$ does not vanish at any point
of $D$, so $\overline{G_1}$ is invertible in $\mathcal{A}$ as claimed.  For
the third assertion, the first part of equation \eqref{dFdG} serves to
define a generator $\overline{\omega}$ in light of parts 1 and 2 above; the
second part of \eqref{dFdG} follows upon expanding the equation
$d(FH_1 + GG_1) = 0$, reducing modulo $\langle F, G \rangle$, and
cancelling $\overline{G_1}$.
\ENDproof

The upshot of the above discussion is that we can represent an element
$\mathcal{A}'$, of the form $\overline{u}\overline{\omega}$ with a unique
$\overline{u} \in \mathcal{A}$, by the column vector $B_u \in \field^3$. In
particular, we represent $\overline{dF} = \overline{G_1\omega}$ by
$B_{G_{1}}$, and $\overline{dG} = \overline{-H_1\omega}$ by
$B_{-H_{1}}$. Hence, we can take the columns of our matrix $M$ to be 
\[
B_{G_1}, B_{xG_1}, B_{yG_1}, B_{-H_1}, B_{-xH_1}.
\]

\begin{proposition}
\label{double-cost2}
Given $F, G, a^{-1}$, the entries of the matrix $M$ can be computed at a
cost of $34M$. 
\end{proposition}
\begin{proof}
We first compute $G_1$ and $H$ at a cost of $10M$, by part 2 of Proposition
\ref{cost7} (recall that we replace $\{F'',G'',G''',H\}$ there by
$\{F,G,G_1,H_1\}$). For later use, we also compute the matrix $T_y$ as in
part 3 of Proposition \ref{cost1}.  This costs us only a further $5M$,
since we have already computed the expression 
$a^{-1}\left (c + d(d-b)\right)$ as part of computing $G_1,H_1$ (when we
computed $(a'')^{-1}\ell$ in the context of the proof of Proposition
\ref{cost7}). As a result, we now have $g$, $h$, and $i$. 

Our next step is to reduce $G_1$ and $H_1$ modulo $\langle F, G \rangle$,
so as to obtain $B_{G_{1}}$ and $B_{H_{1}}$; the extra negation to get
$B_{-H_1}$ costs nothing.  We reduce $G_1 \equiv G_1 - G$ at no
multiplicative cost, and since $G_1 - G \in \field \cdot 1 + \field \cdot x
+ \field \cdot y$ from our formulae for $G_1$ and $G$, we obtain $B_{G_1}$
for free.  As for $H_1$, we have
$H_1 = -y^2 + ax^2 + (\field\text{-linear combination of } 1, x, y)$; hence
by \eqref{s3}
\[
H_1 \equiv H_1 + y^2 + gy + hx + i - aF
             \in \field \cdot 1 + \field \cdot x + \field \cdot y
\]
will be reduced.  The only multiplication needed is to obtain $aF$, which
costs $2M$ to obtain $a^2, ac$, since we already found $ab$ as part of
finding $G_1,H_1$.

Finally, we multiply $T_x$ by the $3\times2$ matrix 
$\left ( B_{G_1} \mid B_{-H_1} \right)$ to obtain $B_{xG_1}$ and
$B_{-xH_1}$ at a cost of $11M$, as in the proof of Proposition
\ref{add-cost2}; we also obtain $B_{yG_1} = T_y B_{G_1}$ at a cost of $6M$,
by part 3 of Proposition \ref{cost1}. 
\ENDproof

\section{Finding the kernel of $M$}
\label{findkernel}

To find $W^{10}_{D+D^{\prime}}$ (respectively $W^{10}_{2D}$) in the case of
addition (respectively, doubling), we must now determine the kernel of our
$3 \times 5$ matrix $M$ from Section \ref{pre-add} (respectively, Section
\ref{pre-double}). A vector
\[
v = \left (
\begin{array}{l}
c_1\\
\vdots\\
c_5
\end{array}
\right)
\] 
satisfying $Mv = 0$ corresponds in both cases to 
\[
L  = c_1 F^{\prime} + c_2 xF^{\prime} + c_3 yF^{\prime}
        + c_4 G^{\prime} + c_5 xG^{\prime} \in W^{10}_{D+D^{\prime}},
\]
since $D = D'$ in the case of doubling.  Our later calculations will be
significantly simplified if we can find a basis $\{s, t\}$ for
$W^{10}_{D+D'}$ of the following special ``monic'' form: 
\[
\begin{array}{rcl}
s & = & x^3 + (\field\mbox{-linear combination of }y^2,xy,x^2,y,x,1)\\
& = & 0x^2y + 1x^3 + \ldots \quad \in W^{9}_{D+D^{\prime}}\\
t & = & x^2y + (\field\mbox{-linear combination of }y^2,xy,x^2,y,x,1)\\
& = & 1x^2y + 0x^3 + \ldots  \quad \in W^{10}_{D+D^{\prime}}.
\end{array}
\]
To do this, we actually find the kernel of a modification $M^{\prime}$ of
$M$: if $M$ has columns  
\[
\left (
\begin{array}{c|c|c|c|c}
 & & & & \\
K_1 & K_2 & K_3  & K_4  & K_5 \\
& & & &  
\end{array}
\right),
\]
then $M^{\prime}$ has columns
\[
\left (
\begin{array}{c|c|c|c|c}
 & & & & \\
K_1 & K_4 & K_3-K_5  & K_2  & K_5 \\
& & & &  
\end{array}
\right).
\]
Note that $M^{\prime}$ can be calculated from $M$ without any field
multiplications. In the case of addition, the columns of $M^{\prime}$
correspond to
\[
\left (
\begin{array}{c|c|c|c|c}
 & & & & \\
\overline{F^{\prime}} & \overline{G^{\prime}}
     & \overline{yF^{\prime}-xG^{\prime}} &
\overline{xF^{\prime}}  & \overline{xG^{\prime}} \\
& & & &  
\end{array}
\right),
\]
and a vector $\left (c_1^{\prime},\ldots, c_5^{\prime}\right)^{{\bf T}} \in
\mbox{ker }M^{\prime}$ corresponds to a combination  
\[
c_1^{\prime}F^{\prime} + c_2^{\prime}G^{\prime}
       + c_3^{\prime}(yF^{\prime}-xG^{\prime})
 + c_4^{\prime}(xF^{\prime})+ c_5^{\prime}(xG^{\prime})
    \in W^{10}_{D+D^{\prime}};
\]
an analogous statement holds in the case of doubling.

We shall see in Section \ref{findST} that the ``monic'' element $s$
comes from a kernel vector with $c_5^{\prime} = 0, c_4^{\prime} = 1$, while
$t$ comes from a kernel vector with $c_5^{\prime} = 1, c_4^{\prime} = 0$.
We thus perform row reduction on $M^{\prime}$ so as to express the unknown
cofficients $c_1^{\prime},c_2^{\prime},c_3^{\prime}$ in terms of the ``free
variables'' $c_4^{\prime}$ and $c_5^{\prime}$.

We write the entries of the modified matrix $M'$ as:
\[
M^{\prime} = \left (
\begin{array}{lllll}
A_1 & B_1 & C_1 & D_1 & E_1\\
A_2 & B_2 & C_2 & D_2 & E_2\\
A_3 & B_3 & C_3 & D_3 & E_3
\label{final3}
\end{array}
\right), 
\]
with rows
$R_i = (A_i \text{ } B_i \text{ } C_i \text{ } D_i \text{ } E_i)$, 
$i = 1, 2, 3$.   
\begin{proposition}
A basis for the kernel of $M^{\prime}$ can be obtained using $39M, 1I$.
\label{cost3}
\end{proposition}
\begin{proof}
Apply row operations to the rows $R_1, R_2, R_3$. This transforms $M'$ into
the following echelon form with the same kernel:
\[
\left (
\begin{array}{lllll}
A_1 & B_1 & C_1 & D_1 & E_1\\
0 & D  & \sigma_1 & \sigma_2 & \sigma_3\\
0 & 0 &  U & \sigma_4  & \sigma_5
\label{final4}
\end{array}
\right ),
\]
where the new rows are 
$R'_1 = R_1,
 R'_2 = A_1 R_2 - A_2 R_1,
 R'_3 = \Delta_{12} R_3 - \Delta_{13} R_2 + \Delta_{23} R_1$.
Here, the $\Delta_{ij}$'s are $2\times2$ minors coming from the first two
columns of $M'$, as given by the formulae below.  This requires us to
compute the following quantities at a cost of $21M$:
\[
\begin{array}{rcl}
D = \Delta_{12} & = & A_1B_2 - A_2B_1\\
\Delta_{13} & = & A_1B_3 - A_3B_1\\
\Delta_{23} & = & A_2B_3 - A_3B_2\\
\sigma_1 & = & A_1C_2 - A_2C_1\\
\sigma_2 & = & A_1 D_2 - A_2 D_1 \\
\sigma_3 & = & A_1 E_2 - A_2 E_1\\
U & = &  \Delta_{12} C_3 - \Delta_{13} C_2 + \Delta_{23} C_1\\
\sigma_4 & = &  \Delta_{12} D_3 - \Delta_{13} D_2 + \Delta_{23} D_1\\
\sigma_5 & = &  \Delta_{12} E_3 - \Delta_{13} E_2 + \Delta_{23} E_1.
\end{array}
\]

\noindent To perform back substitution, we need to obtain 
\begin{equation}
A_1^{-1}, D^{-1}, \quad \mbox{ and } U^{-1}.
\label{inverses}
\end{equation}
For this, we perform
\[
Q_1 = A_1 D, \quad Q_2  =  Q_1 U, \quad Q_3  =  Q_2^{-1},\
\]
\[
U^{-1} = Q_1 Q_3, \quad     Q_4 = U Q_3, \quad
D^{-1} = A_1 Q_4, \quad A_1^{-1}= D Q_4, 
\]
so the inverses in (\ref{inverses}) above can all be produced using
$6M, 1I$.  Back substitution performed on the matrix in (\ref{final4}) now
costs a further $6M + 6M = 12M$ to find the two basis elements
$(\alpha,\beta,\gamma,1,0)^{{\bf T}}$ and 
$(\delta,\varepsilon, \zeta,0,1)^{{\bf T}}$
of the kernel, corresponding to $s$ and $t$. (Solve for $\gamma$, $\beta$,
$\alpha$, $\zeta$, $\varepsilon$, $\delta$ in that order). 
\ENDproof

\section{Finding $s$ and $t$}
\label{findST}

At this point, we have obtained a basis $\{v_1', v_2'\}$ for the kernel of
$M^{\prime}$ of the form
\[
v_{1}^{\prime} = 
\left (
\begin{array}{c}
\alpha\\
\beta\\
\gamma\\
1\\
0
\end{array}
\right),
\]
corresponding to $s$, and 
\[
v_{2}^{\prime} = 
\left (
\begin{array}{c}
\delta\\
\varepsilon\\
\zeta\\
0\\
1
\end{array}
\right),
\]
corresponding to $t$. The desired elements $s$ and $t$ are
\[
\left \{ 
\begin{array}{rcl}
s & = & \alpha F^{\prime} + \beta
G^{\prime} + \gamma(yF^{\prime}-xG^{\prime}) + xF^{\prime}\\
t & = & \delta F^{\prime} + \varepsilon G^{\prime}
            + \zeta (yF^{\prime}-xG^{\prime}) + xG^{\prime}.
\end{array}
\right.
\]
(This includes the case of doubling, for which $F'=F$ and $G'=G$.)
We now have the following:

\begin{proposition}
Given $v_1^{\prime}$ and $v_2^{\prime}$ as above, $s$ and $t$ can be
obtained at a cost of $18M$.
\label{cost4}
\end{proposition}
\begin{proof}
To calculate $s$ and $t$ using as few multiplications as possible, we
illustrate the following steps for $s$ (those for $t$ follow similarly). We
have 
\[
s = (\alpha + \gamma y)F^{\prime} + (\beta - \gamma x)G^{\prime} + xF^{\prime},
\]
where
\[
\begin{array}{rcl}
F^{\prime} & = & x^2 + a^{\prime} y + b^{\prime}x + c^{\prime} \\
G^{\prime} & = & xy + d^{\prime}y + e^{\prime}x + f^{\prime}.
\end{array}
\]
We now wish to expand $s$ as a linear combination of the monomials
$x^3$, $y^2$, $xy$, $x^2$, $y$, $x$, and $1$. Write 
\[
\begin{array}{rcll}
s & = & (\alpha + \gamma y)x^2 + (\beta - \gamma x) xy + xF^{\prime}
& \hskip2em \mbox{(I)} \\
& & + (\alpha + \gamma y)(a^{\prime}y + b^{\prime}x + c)
    + (\beta - \gamma x)(d^{\prime}y + e^{\prime}x + f^{\prime}).
& \hskip2em \mbox{(II)}
\end{array}
\]
The terms in (I) do not involve any multiplication in $K$ (note that the
leading coefficient $x^3$ comes from $xF^{\prime}$). The terms in (II) can
be written as 
\[
\begin{array}{lcll}
(\alpha + \gamma y)b^{\prime}x + (\beta - \gamma x)d^{\prime}y & &
& \hskip9em \mbox{(III)}\\
+ (\alpha + \gamma y )(a^{\prime}y + c)
        + (\beta - \gamma x)(e^{\prime}x + f^{\prime}), & &
&\hskip9em \mbox{(IV)}
\end{array}
\]
where (III) requires $3M$ to form $\gamma(b^{\prime}-d^{\prime})xy
+ \alpha b^{\prime}x + \beta d^{\prime}y$ and (IV) requires $6M$ in total,
using Karatsuba's method for each of the two terms. The total cost is thus
$9M$ to find $s$.

Finding $t$ also requires $9M$; the only essential difference is that
$xF^{\prime}$ becomes $xG^{\prime}$ in the analogue of (I).

The total cost to find $s$ and $t$ is thus $18M$.  Note from the
computation that $s$ and $t$ are both monic in the sense that their
``leading'' coefficient is $1$, and that moreover the coefficient of $x^3$
in $t$ is zero.
\ENDproof

\section{Calculating $xt, yt, x^2t, xyt$ and $xs, ys, x^2s, xys, y^2s$}
\label{xyst}

We have now computed $s,t \in W^{10}_{D+D^{\prime}}$. 
We let $s_1,\ldots,s_6, t_1,\ldots,t_6$ be the coefficients of $s$
and~$t$, as in equation~\eqref{dstar0} above.
As we saw in Subsection
\ref{step2}, we now wish to find $F^{\dprime}, G^{\dprime} \in
W^{7}_{D^{\dprime}}$ via 
\[
\field F^{\dprime} + \field G^{\dprime}
   = \{ \ell \in W^7 \mid \ell t \in sW^8 + W^9\}.
\]
Thus, $\ell$ is a $\field$-linear combination of the basis
$\{t,xt,yt,x^2t,xyt\}$ for $tW^7$ that is congruent to a $\field$-linear
combination of the basis $\{s,xs,ys,x^2s,xys,y^2s\}$ for $sW^8$ in the
quotient space $W^{17}/W^9$. We express these multiples of $s$ and $t$ in
terms of the following ordered basis for $W^{17}$:
\begin{equation}
\{1,x,y,x^2,xy,y^2,x^3,x^2y,xy^2,y^3,x^3y,x^2y^2,xy^3,y^4,x^3y^2\}.
\label{basis17}
\end{equation}
To work in $W^{17}/W^9$, we need only the coefficients of the
last eight monomials: 
\begin{equation}
\{x^2y,xy^2,y^3,x^3y,x^2y^2,xy^3,y^4,x^3y^2\}
\label{smallbasis}
\end{equation}
\begin{lemma}
Given $s$ and $t$ as above, producing the relevant coefficients of $xt$,
$yt$, $x^2t$, $xyt$, $xs$, $ys$, $x^2s$, $xys$, and $y^2s$ requires $2M$.
\label{prelim-cost5}
\end{lemma}
\begin{proof}
Our choice of basis for $W^{17}$ means that we use the equation of the
curve (\ref{eqn}) to eliminate all monomials $x^i y^j$ with $i \geq 4$.
Carrying this out for the multiples of $s$ and $t$ above, we
obtain the matrix
\begin{equation*}
N = \left( 
\begin{array}{ccccccccccc}
t_6 & 0 & 0 & t_3q_0 & 0 & s_6 & q_0 & 0 & s_3q_0 & 0 & 0 \\ 
t_5 & t_6 & 0 & t_3q_1 & 0 & s_5 & s_6+q_1 & 0 & q_0+s_3q_1 & 0 & 0 \\ 
t_4 & 0 & t_6 & q_0+t_3p_0 & 0 & s_4 & p_0 & s_6 & s_3p_0 & q_0 & 0 \\ 
t_3 & t_5 & 0 & t_6+t_3q_2 & 0 & s_3 & s_5+q_2 & 0 & s_6+q_1+ s_3q_2 & 0 & 0 \\ 
t_2 & t_4 & t_5 & q_1+t_3p_1 & t_6 & s_2 & s_4+p_1 & s_5 & p_0+s_3p_1 & s_6 + q_1 & 0 \\ 
t_1 & 0 & t_4 & p_0 & 0 & s_1 & 0 & s_4 & 0 & p_0 & s_6 \\ 
0 & t_3 & 0 & t_5 & 0 & 1 & s_3 & 0 & s_5+q_2 & 0 & 0 \\ 
1 & t_2 & t_3 & t_4+q_2+t_3p_2 & t_5 & 0 & s_2+p_2 & s_3 & s_4+p_1+s_3p_2 & s_5 + q_2 & 0 \\ 
0 & t_1 & t_2 & p_1 & t_4 & 0 & s_1 & s_2 & 0 & s_4+p_1 & s_5 \\ 
0 & 0 & t_1 & t_3 & 0 & 0 & 1 & s_1 & s_3 & 0 & s_4 \\ 
0 & 1 & 0 & t_2 & t_3 & 0 & 0 & 1 & s_2+p_2 & s_3 & 0 \\ 
0 & 0 & 1 & t_1+p_2 & t_2 & 0 & 0 & 0 & s_1 & s_2 + p_2 & s_3 \\ 
0 & 0 & 0 & 0 & t_1 & 0 & 0 & 0 & 1 & s_1 & s_2 \\ 
0 & 0 & 0 & 1 & 0 & 0 & 0 & 0 & 0 & 1 & s_1 \\ 
0 & 0 & 0 & 0 & 1 & 0 & 0 & 0 & 0 & 0 & 1 
\label{original}
\end{array}
\right)
\end{equation*}
whose columns represent in that order $t$, $xt$, $yt$, $x^2t$, $xyt$, $s$,
$xs$, $ys$, $x^2s$, $xys$, and $y^2s$ with respect to our full basis for
$W^{17}$ given in (\ref{basis17}) above. However, since we only need the
last eight rows of $N$ to indicate the values in $W^{17}/W^9$, we only need
to work with the matrix
\begin{equation*}
N' = \left( 
\begin{array}{ccccccccccc}
1 & t_2 & t_3 & t_4+q_2+t_3p_2 & t_5 & 0 & s_2+p_2 & s_3 & s_4+p_1+s_3p_2 & s_5 + q_2 & 0 \\ 
0 & t_1 & t_2 & p_1 & t_4 & 0 & s_1 & s_2 & 0 & s_4+ p_1 & s_5 \\ 
0 & 0 & t_1 & t_3 & 0 & 0 & 1 & s_1 & s_3 & 0 & s_4 \\ 
0 & 1 & 0 & t_2 & t_3 & 0 & 0 & 1 & s_2+p_2 & s_3 & 0 \\ 
0 & 0 & 1 & t_1+p_2 & t_2 & 0 & 0 & 0 & s_1 & s_2 + p_2 & s_3 \\ 
0 & 0 & 0 & 0 & t_1 & 0 & 0 & 0 & 1 & s_1 & s_2 \\ 
0 & 0 & 0 & 1 & 0 & 0 & 0 & 0 & 0 & 1 & s_1 \\ 
0 & 0 & 0 & 0 & 1 & 0 & 0 & 0 & 0 & 0 & 1 
\end{array}
\right).
\end{equation*}
This shows that we only need to compute the multiples $t_3\cdot p_2$ and
$s_3\cdot p_2$, thereby proving our result.
\ENDproof

\section{Finding $F^{\dprime},G^{\dprime}$ that span the subspace
$W^{7}_{D^{\dprime}}$}
\label{last}

We refer to the columns of $N'$ above as
\[
N^{\prime} = \left (
      \begin{array}{c|c|c|c|c}
   C_1 & C_2 & C_3 & \ldots & C_{11}
      \end{array} \right).
\]
We now need to find a linear combination of the first five columns $C_1,
\dots, C_5$ of $N^{\prime}$, corresponding to a basis for the image of
$tW^{7}$ in $W^{17}/W^9$, which belongs to the span of the last six columns
$C_6, \dots, C_{11}$ of $N^{\prime}$, corresponding to the image of
$sW^{8}$ in $W^{17}/W^9$.  Let $V$ denote the $5$-dimensional subspace of
$\field^8$ spanned by the columns $C_6,\ldots,C_{11}$ (of course the zero
column $C_6$ is irrelevant), and let $\mathcal{T}$ denote the set of
columns $\{C_1,\ldots,C_5\}$: we thus want to find combinations of columns
of $\mathcal{T}$ that map to zero in the $3$-dimensional quotient
$\field^8/V$. This quotient can be identified with the subspace
$V' \subset \field^8$ given by 
\[
V^{\prime} = \left \{
   \left ( \alpha , \beta , 0 , 0 , \gamma , 0 , 0 , 0  \right)^{{\bf T}}
 \mid  \alpha, \beta, \gamma \in \field
\right\},
\]
since $V$ and $V^{\prime}$ are complementary subspaces. Our first goal is
then to reduce the columns of $\mathcal{T}$ modulo $V$, so as to obtain
elements $\overline{C_1},\ldots,\overline{C_5} \in V^{\prime}$ with
\[
C_i \equiv \overline{C_i} \bmod V 
     = \left ( \begin{array}{c}
\alpha_{i} \\ \beta_{i} \\ 0 \\ 0 \\ \gamma_{i} \\ 0 \\ 0 \\ 0
               \end{array}
        \right), \quad i = 1, \ldots, 5.
\]
After that, we will need to determine the kernel of the $3 \times 5$ matrix 
\[
M^{\dprime}   =  \left ( \begin{array}{ccccc} 
\alpha_1 & \alpha _2 & \alpha_3 & \alpha_4 & \alpha _5\\
\beta_1 & \beta_2 & \beta_3 & \beta_4 & \beta_5\\
\gamma_1 & \gamma_2  & \gamma_3 & \gamma_4 & \gamma_5
\end{array}
\right)
\]
to obtain $F^{\dprime}$ and $G^{\dprime}$. 
\begin{lemma}
Given the matrix $N'$, the columns of $\mathcal{T}$ can be reduced modulo
$V$ to produce the columns of the matrix $M''$,
at a total cost of $19M$. 
\label{prelim-reductions}
\end{lemma}
\begin{proof}
As a preliminary calculation, we find elements $D_8$, $D_{10}$, and
$D_{11}$ of $V$, corresponding respectively to $ys - s_1xs = (y-s_1x)s$,
$x(y-s_1x)s$, and $y(y - s_1x)s$.  This will aid us in reducing columns of
$\mathcal{T}$ modulo $V$.  We have:
\[
D_8 = C_8 - s_1C_7 
    = \left (\begin{array}{c}
              s_3-s_1(s_2+p_2)\\s_2-s_1^2\\0\\1\\0\\0\\0\\0
      \end{array}\right),
\]
\[
D_{10} = C_{10}-s_1C_9
       = \left (\begin{array}{c}
                  s_5 + q_2 - s_1(s_4+p_1+s_3p_2)\\s_4 + p_1\\-s_1s_3\\
                  s_3-s_1(s_2+p_2)\\s_2+p_2-s_1^2\\0\\1\\0
         \end{array}\right),
\]
\[
D_{11} = C_{11}-s_1C_{10}
       = \left (\begin{array}{c}
                  -s_1(s_5+q_2)\\s_5-s_1(s_4+p_1)\\s_4\\-s_1s_3\\
                  s_3-s_1(s_2+p_2)\\s_2-s_1^2\\0\\1
         \end{array}\right).
\]
Calculating $D_8$, $D_{10}$ and $D_{11}$ costs $6M$, as we already know
$s_3p_2$ from $N'$, so it suffices to calculate
\[
s_1(s_2+p_2), \quad 
s_1^{2}, \quad
s_1(s_4+p_1+s_3p_2), \quad
-s_1s_3, \quad
-s_1(s_5+q_2), \quad
{-s_1}(s_4+p_1).
\]
It is clear that $V$ is spanned by $\{C_7,D_8,C_9,D_{10},D_{11}\}$. We now
compute the reduction of columns of $\mathcal{T}$ modulo $V$. 

\noindent First, note that
\[
\overline{C_1} = C_1 \in V^{\prime}
\]
which comes at no cost, so we obtain
\[
\left ( \begin{array}{c}
   \alpha_1 \\ \beta_1 \\ \gamma_1
        \end{array}
\right)
 = \left ( \begin{array}{c}
      1 \\ 0 \\ 0
           \end{array}
 \right). 
\]
Second, 
\[
\overline{C_2} = C_2 - D_8 \in V^{\prime}
\]
which also comes at no cost, so that 
\[
\left ( \begin{array}{c}
   \alpha_2 \\ \beta_2 \\ \gamma_2
        \end{array}
\right)
 = \left ( \begin{array}{c}
       t_2-s_3+s_1(s_2+p_2) \\ t_1-s_2+s_1^2 \\ 0
           \end{array}
    \right). 
\]
Third, we have
\[
\overline{C_3} = C_3 - t_1C_7 \in V^{\prime}, 
\]
costing $2M$ to calculate $t_1C_7$, and hence
\[
\left ( \begin{array}{c}
   \alpha_3 \\ \beta_3 \\ \gamma_3
        \end{array}
\right)
 = \left ( \begin{array}{c}
      t_3-t_1(s_2+p_2)\\ t_2-t_1s_1\\1
            \end{array}
   \right). 
\]
Fourth and fifth, note that 
\[
C_4 - D_{10} = \left ( \begin{array}{c}
                     m_1\\m_2\\m_3\\m_4\\m_5\\0\\0\\0
                        \end{array}
               \right),
  \text{ with } m_i \in \field, \quad i = 1,\ldots,5,
\]
\[
C_5 - D_{11} = \left ( \begin{array}{c}
                     z_1\\z_2\\z_3\\z_4\\z_5\\z_6\\0\\0
                        \end{array}
               \right),
   \text{ with } z_i \in \field, \quad i = 1,\ldots,6,
\]
so that
\[
C_5-D_{11}-z_6 C_9 = \left ( \begin{array}{c}
                     \ell_1\\ \ell_2\\ \ell_3\\ \ell_4\\ \ell_5\\ 0\\ 0\\ 0
                             \end{array}
                    \right),
    \text{ with } \ell_i \in \field, \quad i = 1,\ldots,6.
\]
Hence our desired reductions are
\[
\begin{split}
\overline{C_4} & =  C_4 - D_{10} - m_4 D_8 - m_3 C_7\\
\overline{C_5} & = C_{5}-D_{11} - z_6 C_9 - \ell_4 D_8 - \ell_3C_7,
\end{split}
\]
which ensures that $\overline{C_4}$ and $\overline{C_5}$ belong to $V'$.
We require $4M$ to find $z_6 C_9$, which allows us to calculate the
vectors $C_4 - D_{10}$ and $C_5-D_{11}-z_6 C_9$.  The expressions
$m_4 D_8 + m_3 C_7$ and  $\ell_4 D_8 + \ell_3 C_7$ can now be obtained
simultaneously as the matrix product
\[
\left ( \begin{array}{c|c} C_7 & D_8 \end{array} \right)
\left ( \begin{array}{cc} m_3 & \ell_3 \\ m_4 & \ell_4 \end{array} \right).
\]
The entries of $C_7$ and $D_8$ are mostly zeros and ones, and the only part
of the above matrix product that involves nontrivial multiplications in
$\field$ is the top $2\times2$ submatrix multiplication
\begin{equation}
\left ( \begin{array}{cc}
    s_2+p_2 & s_3-s_1(s_2+p_2)  \\  s_1 & s_2-s_1^2
        \end{array}
 \right)
\left ( \begin{array}{cc} m_3 & \ell_3 \\ m_4 & \ell_4\end{array}\right).
\label{matrixproduct2x2}
\end{equation}
This costs $7M$ using Strassen's technique. At this point we need no
further multiplications to produce $\overline{C_4}$ and $\overline{C_5}$.

Adding up the costs to produce all of $\overline{C_1},\ldots,\overline{C_5}$
concludes the proof.
\ENDproof
\begin{lemma}
Given $s$ and $t$, the columns of $\mathcal{T}$ can be obtained and reduced
modulo $V$, thereby obtaining the matrix $M''$,
at a total cost of $20M$.  (I.e., we can save one multiplication
compared to using Lemmas \ref{prelim-cost5} and~\ref{prelim-reductions}.) 
\label{reductions}
\end{lemma}
\begin{proof}
We claim that the two multiplications $t_3 p_2$ from Lemma
\ref{prelim-cost5} and $s_1(s_4 + p_1 + s_3 p_2)$ from Lemma
\ref{prelim-reductions} can be replaced with a single multiplication.  To
see this, observe that these two multiplications are used only when we
calculate the first coefficient $m_1$ in 
the column vector
$C_4 - D_{10} = (m_1, m_2, m_3, m_4, m_5, 0, 0, 0)^{\mathbf{T}}$.
Now rearrange  
\[
\begin{split}
m_1 &= t_4 + q_2 + t_3 p_2 - s_5 - q_2 + s_1(s_4 + p_1 + s_3 p_2) \\
    & = t_4 - s_5 + s_1(s_4 + p_1) + (t_3 + s_1 s_3)p_2.
  \end{split}
\]
Since we have already computed $ s_1(s_4 + p_1)$ and $s_1 s_3$ during
Lemma \ref{prelim-reductions}, we see that we can replace the two
multiplications $t_3 p_2$ and $s_1(s_4 + p_1 + s_3 p_2)$ by the single
multiplication $(t_3 + s_1 s_3)p_2$.  This concludes our proof. 
\ENDproof

The following proposition now allows us to find the desired polynomials
\[
\begin{array}{rcl}
F^{\dprime} & = & x^2 + a^{\dprime}y + b^{\dprime} x + c^{\dprime}\\
G^{\dprime} & = & xy + d^{\dprime}y + e^{\dprime} x + f^{\dprime}.
\end{array}
\]

\begin{proposition}
Given $s$ and $t$, the polynomials $F^{\dprime}$ and $G^{\dprime}$, as well
as the inverse $(a'')^{-1}$ , can be obtained using $31M, 1I$.
\label{cost6}
\end{proposition}
\begin{proof}
Recall that the columns of $M^{\dprime}$ represent the reductions of each
of $t$, $xt$, $yt$, $x^2t$ and $xyt$ modulo the multiples of $s$ via the
``reduction modulo $V$'' described above. Hence, by Lemma \ref{reductions},
the matrix $M^{\dprime}$ can be obtained using $20M$, and has the form 
\[
M^{\dprime} = \left ( 
\begin{array}{ccccc}
1 & \alpha_2 & \alpha_3 & \alpha_4 & \alpha_5\\
0 & \beta_2 & \beta_3 & \beta_4 & \beta_5\\
0 & 0 & 1 & \gamma_4 & \gamma_5
\end{array}
\right).
\] 
In anticipation of our next step, we compute $\gamma_4^{-1}$ and
$\beta_2^{-1}$ using $3M, 1I$ (i.e., find $\beta_2 \cdot \gamma_4$, invert
it, and multiply the inverse separately with each of  $\beta_2$ and
$\gamma_4$). 
We now can find two vectors 
\[
\left \{
\begin{array}{rcl}
v_1^{\dprime} & = & \left (
        c^{\dprime}, b^{\dprime}, a^{\dprime}, 1,0
                    \right)^{{\bf T}}\\
v_2^{\dprime} & = & \left (
        f^{\dprime}, e^{\dprime}, d^{\dprime},0,1
                    \right)^{{\bf T}}
\end{array}
\right.
\]
that span the kernel of $M^{\dprime}$ using back substitution, requiring a
further $8M$. Those give us the coefficients of the polynomials $F''$ and
$G''$.  Note that $a'' = - \gamma_4$, and so we know its inverse thanks to
our previous anticipatory step.
\ENDproof

\section{Negating the final result, and an application to Section \ref{pre-double}}
\label{negation}

As mentioned in Subsection \ref{step3}, our final result representing
$\xi''' = -\xi'' = \xi + \xi'$ will be a pair $\{F''',G'''\}$,  with
$F''' = F''$ and $G''' = xy + d'''y + e''' x + f''' \in W^7_{D'''}$ that
satisfies $G''G'''+F''H=0$ for some $H \in W^8$.  We can then in principle
find $G'''$ by a procedure analogous to that in Sections \ref{xyst} and
\ref{last}, by working modulo $W^6$, which is analogous to how we
previously dropped some rows from the matrix $N$ to get $N'$.  If we 
furthermore need to find $H$, as is the case in Proposition
\ref{double-cost2}, we can do something similar by dropping one fewer row
at the start, i.e., by working modulo $W^4$ (we invite the reader to check
that this extra ``precision'' is required exactly to obtain the constant
term of $H$). 

We however preferred to find the following solution by a direct calculation:
\begin{equation}
\label{GnegH}
\begin{array}{rcl}
G''' &=& xy + (b''-d'')y - (\ell (a'')^{-1} + m)x \\
      & & \qquad
             + [md'' + (\ell (a'')^{-1}+e'')(d''-b'')
                     + a''(a''b'' - p_1) - f''
               ] \\
H &=& -y^2 + a''x^2 + \ell (a'')^{-1}y - a''b''x \\
      & & \qquad
            +[(\ell (a'')^{-1} + m)e'' + a''({b''}^2 - c'' - q_2)],
\end{array}
\end{equation}
where
 \[
\begin{array}{rcl}
m & = & e'' + a''(a'' + p_2)\\
\ell & = & c'' + (d''-b'')d''.
\end{array}
\]
This can be verified without setting up a system of linear equations:
instead, note that our expressions for $G''', H$ satisfy $G''G''' + F''H
\in \field\cdot x + \field\cdot y + \field\cdot 1 = W^4$ (taking into
account equation (\ref{eqn}) of our curve).  However, any combination of
$F''$ and $G''$ vanishes at $D''$, so $G''G''' + F''H \in W^4_{D''} = 0$,
since $D''$ is typical. 

Thus our result is:
\begin{proposition}
\label{cost7}
Given $F'', G''$ and $(a'')^{-1}$, let $F''', G'''$ represent the negative
in the Jacobian; then $F''' = F''$, and we obtain
$(a'')^{-1} = (a''')^{-1}$ for free.   
\begin{enumerate}
\item
It costs $7M$ to compute $G'''$ as given by the above formulae.
\item
It costs $10M$ to compute both $G'''$ and $H$, satisfying $G''G''' + F''H = 0$.
\end{enumerate}
\end{proposition}
\begin{proof}
First compute $m$ and $\ell$, then compute $\ell (a'')^{-1}$ and $a''b''$,
and then compute the remaining coefficients of $G'''$ (and of $H$, if
needed) using the above expressions.
\ENDproof

\section{Conclusion}
\label{conclusion}

We now assemble all the parts to obtain the main result of our paper:
\begin{theorem}
In the Jacobian of a $C_{3,4}$ curve defined over a large finite field
$\field$, point addition can be performed on typical elements using $117$
field multiplications and $2$ field inversions. Point doubling can be
performed on typical elements using $129$ field multiplications and $2$
field inversions.  
\end{theorem}
\begin{proof}
For point addition, add up the costs of Propositions \ref{add-cost2},
\ref{cost3}, \ref{cost4}, \ref{cost6}, and part 1 of Proposition
\ref{cost7}. For point doubling, add up the costs of Propositions
\ref{double-cost2}, \ref{cost3}, \ref{cost4}, \ref{cost6}, and part 1 of
Proposition \ref{cost7}. 
\ENDproof

In terms of the number of multiplications required, our results represent
improvements of $19.3\%$ for addition and $22.8\%$ for doubling
(compared to \cite{FOR04}), and of $22\%$ for addition and $25.8\%$ for
doubling (compared to \cite{BEFG}). All the algorithms require two
inversions in $\field$ per group operation in the Jacobian. 


\begin{thebibliography}{1}

\bibitem{base-2}
Roberto Avanzi, Gerhard Frey, Tanja Lange, and Roger Oyono.
\newblock On using expansions to the base of {$-2$}.
\newblock {\em Int. J. Comput. Math.}, 81(4):403--406, 2004.

\bibitem{BEFG}
Abdolali Basiri, Andreas Enge, Jean-Charles Faug{\`e}re, and Nicolas G{\"u}rel.
\newblock Implementing the arithmetic of {$C\sb {3,4}$} curves.
\newblock In {\em Algorithmic number theory}, volume 3076 of {\em Lecture Notes
  in Comput. Sci.}, pages 87--101. Springer, Berlin, 2004.

\bibitem{BEFG-article}
Abdolali Basiri, Andreas Enge, Jean-Charles Faug{\`e}re, and Nicolas G{\"u}rel.
\newblock The arithmetic of {J}acobian groups of superelliptic cubics.
\newblock {\em Math. Comp.}, 74(249):389--410 (electronic), 2005.

\bibitem{Diem}
Claus Diem.
\newblock An index calculus algorithm for plane curves of small
degree.
\newblock In {\em Algorithmic number theory}, volume 4076 of {\em Lecture Notes
  in Comput. Sci.}, pages 543--557. Springer, Berlin, 2006.

\bibitem{DT06}
Claus Diem and Emmanuel Thom{\'e}.
\newblock Index calculus in class groups of non-hyperelliptic curves of genus
  three.
\newblock 2006 preprint, to appear in {\em J. of Cryptology}, may be
downloaded from the web at the URL
  \texttt{http://www.math.uni-leipzig.de/$\sim$diem/preprints/non-he-genus3.\{dvi,ps\}}

\bibitem{FO-previous}
St{\'e}phane Flon and Roger Oyono.
\newblock Fast arithmetic on {J}acobians of {P}icard curves.
\newblock In {\em Public key cryptography---PKC 2004}, volume 2947 of {\em
  Lecture Notes in Comput. Sci.}, pages 55--68. Springer, Berlin, 2004.

\bibitem{FOR04}
St{\'e}phane Flon, Roger Oyono, and Christophe Ritzenthaler.
\newblock Fast addition on non-hyperelliptic genus 3 curves.
\newblock 2004 preprint, may be electronically downloaded from the web at
  either of the URLs
  \texttt{http://www.exp-math.uni-essen.de/$\sim$oyono/Quartic.html} or
  \texttt{http://www.math.uwaterloo.ca/$\sim$royono/Quartic.html}

\bibitem{Hess}
Florian Hess.
\newblock {\em Zur {D}ivisorenklassengruppenberechnung in globalen
  {F}unktionenk{\"o}rpern}.
\newblock PhD thesis, Technische Universit{\"a}t Berlin, 1999.
\newblock May be electronically downloaded from the web at
  \texttt{http://www.math.tu-berlin.de/$\sim$kant/publications/diss/\{diss\_FH.ps.gz,hess.pdf\}}

\bibitem{KKM}
Kamal Khuri-Makdisi.
\newblock Linear algebra algorithms for divisors on an algebraic curve. 
\newblock {\em Math. Comp.}, 73(245):333--357 (electronic), 2004; 
\newblock \texttt{math.NT/0105182}

\bibitem{KKM04}
Kamal Khuri-Makdisi.
\newblock Asymptotically fast group operations on {J}acobians of
  general curves.
\newblock Preprint, \texttt{math.NT/0409209} to appear in Mathematics
  of Computation.

\end{thebibliography}

\end{document}